\documentclass[onecolumn]{autart}
\usepackage{graphicx}
\input{amssym}

\renewcommand{\theequation}{\arabic{section}.\arabic{equation}}
\begin{document}
\begin{frontmatter}
\title{Variational problems without having any non-trivial~Lie~variational~symmetries} 
\thanks[footnoteinfo]{ Corresponding author: Tel. +9821-73913426.
Fax +9821-77240472.}
\author[]{M. Nadjafikhah\thanksref{footnoteinfo}}\ead{m\_nadjafikhah@iust.ac.ir},
\author[]{S. Dodangeh}\ead{s\_dodangeh@mathdep.iust.ac.ir},
\address{School of Mathematics, Iran University of Science and Technology, Narmak, Tehran 1684613114, Iran.}
%
\begin{abstract}
In this paper we construct variational problems without Lie
non-trivial variational symmetry and solving them using new class
of symmetries ($\mu$-symmetry) which introduced by Guiseppe Gaeta
and Paola Morando (2004). The central object in this paper is
horizontal one-form $\mu$ on first order jet space $J^1 M$.
\end{abstract}
\end{frontmatter}
\section{Introduction}
Hidden symmetries defined as symmetries that are lost (Type I) or
gained (Type II) as the order of an ODE is reduced or as the
number of variable of a PDE is reduced. Hidden symmetries are
difficult to evaluate since there are no general direct method for
determining them. There are several approach that we can use to
investigate hidden symmetries and gain them.
\newline $~~~~$
In 2001, Muriel and Romero introduced $\lambda$-symmetries to
evaluate Type I hidden symmetries of ODEs \cite{[7]}. Guiseppe
Gaeta and Pola Morando expanded this approach to scalar PDEs and
PDEs systems. They constructed equations without Lie point
symmetries too \cite{[11]}. This equations have no obvious  order
reduction (in ODE case) and variable reduction (in scalar PDEs
case and PDEs systems) which can be reduce using $\mu$-symmetries.
\newline $~~~~$
In this paper we construct equations without Lie non-trivial
symmetries using \cite{[11]}. you can assume these equations are
Euler-Lagrangian of some variational problems (with necessary
condition) and construct the variational problems have this
equations as Euler-Lagrangian, using direct method (This is
inverse problem in variational calculus). For such variational
problems we can't solve them using Lie symmetry method (Lie
classical method), so solve them using this new class of
symmetries ($\mu$-symmetries).
\section{$\mu$-symmetry on scalar PDEs and PDEs systems}
The starting point will be a discussion of some of the
foundational results about $\mu$-symmetry. In this section we
recall these results rather briefly. Reader can consult
\cite{[11]} to gain complete information about this symmetries.
\newline $~~~~$
Let $\mu=\lambda_idx_i$ be horizontal one-form on first order jet
space $(J^1M,\pi,M)$ and compatible with contact structure
$\varepsilon$ on $J^kM$ for $k\geq 2$. i.e.
\begin{eqnarray}\label{eq:1}
d\mu\in J(\varepsilon).
\end{eqnarray}
where $J(\varepsilon)$ is Cartan ideal generated by contact
structure.
\begin{thm}{\rm (See \cite{[11]})}
Condition (\ref{eq:1}) is equivalent to
$D_i\lambda_j-D_j\lambda_i=0$. Where $D_i$ is total derivative
w.r.t. $x_i$.
\end{thm}
Let $X:=\xi^i\partial_{x_i}+\phi(x,u)\partial_u$ be a vector field
on total space $M$ and $q=1$. (i.e the number of dependent
variable is one, on the other hand we discus these concepts in
scalar PDEs framework). We define $Y:=X+\Psi_j\partial u_j$ on
$k$-th order jet space $J^k M$ as $\mu$-prolong of $X$ if its
coefficient satisfy the $\mu$-prolongation formula
\begin{eqnarray}\label{eq:3}
\Psi_{J,i}=(D_i+\lambda_i)\Psi_{J}-u_{J,m}(D_i+\lambda_i)\xi^m.
\end{eqnarray}
\begin{rem}
If we set $\mu=0$ in (\ref{eq:3}) then we gain ordinary
prolongation of $X$. so we can assume ordinary prolong as
$0$-prolong in $\mu$-prolong framework.
\end{rem}
We can show connection between ordinary prolong and $\mu$-prolong
in follow theorem.
\begin{thm}{\rm (See \cite{[11]})}
{\it Let $X:=\xi^i\partial_{x_i}+\Phi(x,u)\partial_u$ be a vector
fields on first order jet space $J^1M$ and
$Y=X+\Psi_j\partial_{u_j}$ be $\mu$-prolong of $X$ and
$X^{(k)}=X+\phi_j\partial_{u_j}$ be ordinary prolong of $X$. Then
we have $\Psi_J=\Phi_J+F_J$, where $F_J$ satisfy the recursion
relation (with $F_0=0$):
$F_{J,i}=(D_i+\lambda_i)F_J+\lambda_iD_iQ$; where $Q$ is Lie
characteristic.} \end{thm}

This theorem provide an economics way of computing
$\mu$-prolongation of $X$ if we knew already its ordinary
prolongation.
%
%
\section{Variational problems and Lie standard reduction method}
Variational problem is finding the extremals (maxima and/or
minima) of a functional
\begin{eqnarray}\label{eq:6}
\ell_\alpha(L)=\sum_J(-D)_J\,\partial_{ u_J^\alpha}.\hspace{2cm}
\alpha=1,2,...,q
\end{eqnarray}
over some space of functions $u=f(x)$, $x\in\Omega$. For such
problems we can define Euler-Lagrangian operators as
\begin{eqnarray}\label{eq:7}
E_\alpha=\sum_J(-D)_J\,\partial_{ u^\alpha_J}.\hspace{2cm}
\alpha=1,2,...,q
\end{eqnarray}
Symmetry on variational problems is motivated by following theorem
of the calculus of variational problems.
\begin{thm}{\rm (See \cite{[4],[9],[3]})}
The smooth extremals $u=f(x)$ of variational problem with
Lagrangian $L(x,u^{(n)})$ must be satisfied in the systems of its
related Euler-Lagrange equations.
\begin{eqnarray}\label{eq:8}
E_\alpha(L)=\Sigma_J(-D)_J\frac{\partial L}{\partial u^\alpha_J}
\end{eqnarray}
\end{thm}
Now we describe our approach in Lie classical method to find
extremals of variational problems on open connected domain
$\Omega$ (see \cite{[4],[5],[9],[3]}). In first step we compute
Euler-Lagrangian equations of problems, next we characterize
symmetry group of this equations and solve this equations using
Lie symmetry method, finally we check this solutions in original
problems. Now, what we can do when Euler-Lagrangian equation has
no standard Lie symmetry? In this paper we gain variational
problems without Lie non-trivial variational symmetry and show how
we can solve such problems. Main theorem of our approach is
following:
\begin{thm}{\rm (See \cite{[2],[3]})}
G is variational symmetry of variational problem (\ref{eq:6}) if
and only if it is Lie symmetry of its Euler-Lagrangian equation
(\ref{eq:7}).
\end{thm}

\section{Euler-Lagrangian equations without Lie non-trivial symmetry}
In first step we characterize equations (scalar PDEs) without Lie
symmetry. For this purpose we consider $X$ be vector field on
first order jet space $J^1M$, then determine general scalar PDEs
with no Lie non-trivial symmetry which admit $X$ as
$\mu$-symmetry.

Consider the vector field
$X=x^2\partial_x+t\partial_t+u\partial_u$. We have
$Q=u-x^2u_x-tu_t$. The corresponding coordinates $(y,v)$ and the
parametric coordinate $\sigma$ in $M=(x,t,u)$ can be chosen as
$\sigma=t$, $y=te^{1/x}$ and $v=u/t$. The corresponding inverse
change of variables is $x=-\ln(y\sigma)$, $t=\sigma$, $u=v\sigma$.
Hence, the function $v=v(\sigma,y)$, is $X$-invariant if and only
if $v_\sigma=0$. The partial derivations of $u$ express in the
partial derivatives of $v=v(\sigma,y)$ as
\begin{eqnarray}\label{eq:11}
u_x=\frac{-y^2}{\ln^2 y\sigma}v_y,\qquad u_t=v+\sigma
v_\sigma+yv_y.
\end{eqnarray}
The above can be inverted to give
\begin{eqnarray}\label{eq:12}
v_y=-\frac{(2\ln t+1/x)^2}{t^2e^{2/x}}u_x,\qquad
v_\sigma=\frac{1}{t}\Bigg[u_t-\frac{u}{t}+\frac{(2\ln
t+1/x)^2}{te^{1/x}}u_x\Bigg];
\end{eqnarray}
Similar above we have this expressions for second order
derivatives.
\begin{eqnarray}\label{eq:13}
&&u_{xx}=-y\ln^2
({y/\sigma})v_{yy}+\frac{2y^2\ln^2({y/\sigma})+2\ln(y\sigma)\ln^2
({y/\sigma})}{\ln^4(y\sigma)}v_y,\nonumber\\
&&u_{xt}=-y\ln^2 ({y/\sigma})\hspace{.1cm}v_y-y^2\ln^2(y/\sigma)v_{yy},\\
&&u_{tt}=-\frac{v}{\sigma}+v_\sigma+\sigma
v_{\sigma\sigma}+2yv_{\sigma
y}+\frac{y}{\sigma}v_y+\frac{y^2}{\sigma}v_{yy},\nonumber
\end{eqnarray}
As before object in this computation is horizontal one-form $\mu$
on one order jet space. In this case since independent variables
is two-dimensional as a result we have: $\mu=\lambda dx+\tau dt$.
\newline
Let us come to the second $\mu$-prolongation of $X$. (Standard
prolongation will be ordinary by setting $\lambda=\tau=0$).
\newline
For this computation we can use  (\ref{eq:3}) or recursion
relation in theorem 3. Hence if we show $\mu$-prolongation of $X$
as
\begin{eqnarray}\label{eq:14}
Y=X+\Psi^x\partial_{u_x}+\Psi^t\partial_{u_t}+\Psi^{xx}\partial_{u_{xx}}+\Psi^{xt}\partial_{u_{xt}}+\Psi^{tt}\partial_{u_{tt}}
\end{eqnarray}
Then we have
\begin{eqnarray}\label{eq:15}
&& \Psi^{x}=(1-2x)+\lambda Q,\qquad \Psi^{t}=\tau Q,\nonumber\\
&& \Psi^{xx}=(1-4x)u_{xx}-2u_x+2\lambda(D_xQ)+[\lambda^2+(D_x\lambda)]Q\\
&& \Psi^{xt}=-2xu_{xt}+[\lambda(D_tQ)+\tau(D_xQ)]+(1/2)[2\lambda\tau+(D_t\lambda)+(D_x\tau)]Q,\nonumber\\
&& \Psi^{tt}=-u_{tt}+2\tau(D_tQ)+[\tau^2+(D_t\tau)]Q.\nonumber
\end{eqnarray}
We consider two simplest case for $\mu$ instead of general case.
\newline
\textbf{Case I:}
 $\tau=0$ and $\lambda$ is real number.
\newline
In this case by substituting this $\mu$ in above we find,
\begin{eqnarray}\label{eq:16}
&& \Psi^x=(1-2x)+\lambda(u-x^2u_x-tu_t),\qquad \Psi^t=0,\nonumber \\
&& \Psi^{xx}=(1-4x)u_{xx}-2u_x-2\lambda[(2x-1)u_x+x^2u_{xx}+tu_{xt}]+\lambda^2(u-x^2u_x-tu_t),\\
&& \Psi^{xt}=-2u_{xt}-\lambda(x^2u_{xt}+tu_{tt}),\qquad
\Psi^{tt}=-u_{tt},\nonumber
\end{eqnarray}
Now if we take $(y,v)$ as invariants of order zero, $\xi_1,\xi_2$
invariants of order one and $(\eta_1,\eta_2,\eta_3)$ invariants of
order two then we find
\begin{eqnarray}\label{eq:17}
&& y=te^{1/x},\qquad v=\frac{u}{x},\qquad \xi_1=\ln2-\frac{1}{2}\ln(-x^2\lambda(S_1-uS_2)^2)+\frac{1}{2S_2}(\lambda-S_2)\ln(\frac{S_1+uS_2}{S_1-uS_2})\nonumber \\
&& \xi_2=u_t,\qquad \eta_1=\ln2-\frac{1}{2}\ln(-(S_3+uS_4)^2(1-4x+2\lambda x^2))+\frac{1}{2}(\frac{\lambda^2}{S_4}-1)\ln(-\frac{S_3-uS_4}{S_3+uS_4}),\\
&& \eta_2=\frac{1}{2}\lambda
x^2u_{xt}^2-\frac{1}{2}u^2+u_{xt}^2-tu_{utt}u_{xt},\qquad
\eta_3=tu_{tt},\nonumber
\end{eqnarray}
Where $S_1,S_2,S_3$ and $S_4$ are respectively
\begin{eqnarray}\label{eq:18}
&& S_1=2\lambda x^2u_x-2+4x+2\lambda tu_t+\lambda u,\qquad S_2=\sqrt{\lambda(4x^2+\lambda)},\qquad S_4=\sqrt{-4+16x-8\lambda x^2+\lambda^4}, \\
&& S_3=2u_{xx}-8xu_{xx}+4\lambda x^2u_{xx}+4u_x+8\lambda
xu_x-4\lambda u_x+4\lambda
tu_{xt}+2\lambda^2x^2u_x+2\lambda^2tu_t-\lambda^2u.\nonumber
\end{eqnarray}
\begin{thm}
Consider the equation
$\Delta:=F(y,v,\xi_1,\xi_2,\eta_1,\eta_2,\eta_3)$ with arbitrary
smooth function $F$. Let $\lambda$ be a real constant.
Then\\
i) The equation $\Delta$ admits the vector field $X$ as a
$\mu$-symmetry with $\mu=\lambda dx$. \\
ii) For $(\partial F/\partial\xi_1)^2+(\partial
F/\partial\eta_1)^2+(\partial F/\partial\eta_2)^2\neq0$, $X$ is
not an ordinary symmetry of $\Delta$.
\end{thm}
\begin{pf} {\rm i)}
As mentioned Lie point symmetry method, PDE equation admit $X$ as
$\mu$-symmetry when we can rewrite it in terms of $X$-invariants
\cite{[11]}. Hence equation $\Delta$ admits $X$ as $\mu$-symmetry
with $\mu=\lambda dx$.

{\rm ii)} Using (\ref{eq:17}), we conclude, $\xi_1,\eta_1$ and
$\eta_2$ depend on $\mu$  in solution space ($I_x$). So if $F$
depend on this arguments then $X$ is not ordinary symmetry of
$\Delta$.
\end{pf}
\textbf{Case II:} $\lambda=0$ and $\tau$ is real number
\newline
Now by substituting this equation in (4.9), we have
\begin{eqnarray}\label{eq:19}
&& \Psi^x=(1-2x),\qquad \Psi^t=\tau(u-x^2u_x-tu_t),\qquad \Psi^{xx}=(1-4x)u_{xx}-2u_x,\nonumber \\
&& \Psi^{xt}=-2xu_{xt}+\tau(u_x-2xu_x-x^2u_{xx}-tu_{xt}),\\
&&
\Psi^{tt}=-u_{xt}+2\tau(u_t-x^2u_{xt}-u_t-tu_{tt})+\tau^2(u-x^2u_x-tu_t),\nonumber
\end{eqnarray}
So we find
\begin{eqnarray}
&& y=te^{1/x},\qquad v=\frac{u}{t},\qquad \xi_1=-\frac{1}{2}\frac{t^2-4xu_x+2u_x}{2x-1},\nonumber \\
&& \xi_2=-\frac{1}{2}\ln\Big(-\frac{1}{4}\frac{t(S_1+uS_2)^2}{\tau}\Big)-\frac{1}{2}\ln(-\frac{S_1-uS_2}{S_1+uS_2})-\frac{1}{2}\frac{\tau}{S_2}\ln(\frac{S_1-uS_2}{S_1+uS_2}),\\
&& \eta_1=\frac{1}{2}(4x-1)u_{xx}^2+4u_xu_{xx}-t^2,\nonumber \\
&& \eta_2=\frac{1}{2}\tau tu_{xt}^2-\tau u_xu_{xt}+\tau x^2u_{xx}u_{xt}+2\tau xu_xu_{xt}-\frac{1}{2}u^2+xu_{xt}^2\nonumber \\
&&
\eta_3=\frac{1}{2}(\frac{\tau^2}{S_4}-1)\ln(-\frac{S_3-uS_4}{S_3+uS_4})-\frac{1}{2}\ln(-\frac{1}{2}t\tau(S_3+S_4)^2),\nonumber
\end{eqnarray}
Where $S_1,S_2,S_3$ and $S_4$ are respectively,
\begin{eqnarray}\label{eq:20}
&& S_1=2\tau tu_t+2\tau x^2u_x-\tau u\qquad S_2=\sqrt{\tau(\tau+4t)}\\
&& S_3=-4\tau tu_{tt}+2u_{xt}+4\tau
x^2u_{xt}+2\tau^2tu_t-u\tau^2,\qquad
S_4=\sqrt{\tau(\tau^3+8t)},\nonumber
\end{eqnarray}
Similar preceding theorem, we have
\begin{thm}
Consider the equation
$\Delta:=F(y,v,\xi_1,\xi_2,\eta_1,\eta_2,\eta_3)$ with arbitrary
smooth function $F$. Let $\tau$ be a real constant. Then\\ i) The
equation $\Delta$ admits $X$ as $\mu$-symmetry with $\mu=\lambda
x$.\\ ii) For $(\partial F/\partial\xi_2)^2+(\partial
F/\partial\eta_2)^2+(\partial F/\partial\eta_3)^2\neq0$, $X$ is
not ordinary symmetry of $\Delta$.
\end{thm}
Using such a procedure we can construct scalar PDEs with
$\mu$-symmetries. If we set
$X=\sum_i\xi^i\partial_{x_i}+\partial_u$ and apply the mentioned
procedure, then we get (p+1)-PDEs without Lie non-trivial
symmetries which have $X$ as $\mu$-symmetry. We can solve this
equations similar to Lie standard symmetry method using $X$ as new
symmetry.

Now we express step II for construct our favorite variational
problems.
\section{$\mu$-symmetry on variational problems}
Characterizing systems of differential equations which are the
Euler-Lagrange equations for some variational problems, is known
as the inverse problem in the calculus of variations (see
\cite{[3],[9]}). There are different approaches to solve or
investigate inverse problem (\cite{[2],[9]}). In order to keep the
scope manageable, we use direct method in this paper.

In this section we construct two examples of variational problems
which their Euler-Lagrange equations have no Lie standard
symmetry. For this purpose first we assume $\Delta$ be equations
without Lie symmetry then find appropriate variational problem
which have $\Delta$ as Euler-Lagrange equation.
\begin{exmp}
Consider the equation
\begin{eqnarray}\label{eq:21}
u_{xx}=[(x+x^2)e^u]_x
\end{eqnarray}
This equation appear in page 182 of P.J. Olver \cite{[3]} as an
equation which can be integrated by quadratures, but lacks
non-trivial symmetries.
\newline
Muriel and Romero in \cite{[6]} solve this equation by using
$\lambda$-symmetry with $\lambda=[(x+x^2)e^u]_u$ and $X=\partial_{
u}$. Let this equation be Euler-Lagrange equation of some second
order variational problem, so we have
\begin{eqnarray}\label{eq:22}
\frac{\partial L}{\partial x}+u_x\frac{\partial L}{\partial
u}+u_{xx}\frac{\partial L}{\partial u_x}+u_{xxx}\frac{\partial
L}{\partial u_{xx}}=u_{xx}-(1+2x)e^u-(x+x^2)u_xe^u,
\end{eqnarray}
Where by solving this equation we find following Lagrangian:
\begin{eqnarray}\label{eq:23}
L(x,u,u_x,u_{xx})=-x^2e^u-e^ux+xu_{xx}+F(u_{xx},u_x-xu_{xx},u+\frac{1}{2}x^2u_{xx}-xu_x);
\end{eqnarray}
Where $F$ is an arbitrary function.
\end{exmp}
So we have following proposition using (theorem 5)and and
corollary 7.4 in \cite{[3]}:
\begin{prop}
The following variational problem and any variational problem with
lagrangian $\hat{L}=L+Div\xi$ with arbitrary smooth function $\xi$
have no Lie non-trivial variational symmetry, and its
Euler-Lagrangian equation has $\lambda$-symmetry with
$\lambda=[(x+x^2)e^u]_u$ and $X=\partial_{ u}$.
\begin{eqnarray}\label{eq:24}
\ell(u)=\int_{\Omega}(-x^2e^u-e^ux+xu_{xx}+F(u_{xx},u_x-xu_{xx},u+\frac{1}{2}x^2u_{xx}-xu_x)).dx,
\end{eqnarray}
where $F$ is an arbitrary function.
\end{prop}
\begin{exmp}
Consider this equation
\begin{eqnarray}\label{eq:25}
8(u_x+1)u_{xx}-24(xu_x^2)-2(u^2x^2+2u_x+24u+1)u_x+x^3u^5+(5x^2+8x)u^4+(7x+32).u^3+3u^2=0,
\end{eqnarray}
Muriel and Romero in \cite{[6]} prove that this equation has no
Lie non-trivial symmetry. Now we use direct method to find some
second variational problem with property of equation (5.11) which
is its Euler-Lagrange equation (\ref{eq:25}).
\begin{eqnarray}\label{eq:26}
&&\hspace{-1cm}\frac{\partial L}{\partial x}+u_x\frac{\partial
L}{\partial u}+u_{xx}\frac{\partial L}{\partial
u_x}+u_{xxx}\frac{\partial
L}{\partial u_{xx}}=\nonumber\\
&=&8(u_x+1)u_{xx}-24(xu_x^2)-2(u^2x^2+2u_x+24u+1)u_x+x^3u^5+(5x^2+8x)u^4+(7x+32).u^3+3u^2\\
&=&0,\nonumber
\end{eqnarray}
By solving this equation we have:
\begin{eqnarray}\label{eq:27}
&& \hspace{-5mm}L(x,u,u_x,u_xx)=F(u_x,u-xu_x,u_{xx})+\frac{1}{504}x^9u_x^5+\frac{1}{56}x^8u_x^4u+\Big(\frac{1}{21}u_x^4-\frac{1}{14}u_x^3u^2\Big)x^7\nonumber\\
&& +\Big(-\frac{1}{3}u_x^3u+\frac{4}{15}u_x^4+\frac{1}{6}u_x^2u^3\Big)x^6+\Big(-\frac{5}{12}u_x^3-\frac{8}{5}u_x^3u+u^2u_x^2-\frac{1}{4}u^4u_x\Big)x^5\\
&& +\Big(\frac{1}{4}u^5+\frac{5}{12}u_x^2u-\frac{5}{3}u_xu^3+4u^2u_x^2\Big)x^4+\Big(-\frac{16}{3}u_xu^3+u_x^2-\frac{25}{6}u_xu^2+\frac{5}{3}u^4+32u_x^2u\Big)x^3\nonumber\\
&&
+\Big(4u^4-3u_xu+\frac{7}{2}u^3+12u_x^2-48u_xu^2\Big)x^2+\Big(32u^3+3u^2-48uu_x-4u_x^2+(8u_{xx}-2)u_x+8u_{xx}\Big)x\nonumber
\end{eqnarray}
where $F$ is an arbitrary function.
\end{exmp}
As a result above, (theorem 5) and corollary 7.4 in \cite{[3]}, we
can find this proposition,
\begin{prop}
Variational problem (5.12) and any Variational problem with
$\hat{L}=L+Div\xi$ with arbitrary function $\xi$ have no Lie
nontrivial variational symmetries and its Euler-Lagrangian
equation has $X=u\partial_{u}$ as $\lambda$-symmetry with
$\lambda=x/u^2$.
\end{prop}
\section*{Conclusion}
In this paper first we construct Euler-Lagrange equations with no
Lie non-trivial symmetry, next we find related variational
problems without any Lie non-trivial variational symmetries.
Finally we solve these variational problems using $\mu$-symmetry.


\begin{thebibliography}{9}
\bibitem[Ab-1996]{[5]} {\sc Abraham-Shrauner, B.} {\em Hidden Symmetries and nonlocal group generators for ordinary differential eauations}, IMA J. Appl. Math. 56,
235-252, 1996.
\bibitem[Ge-Mo-2004]{[11]} {\sc Gaeta, G.} and {\sc Morando, P.} {\em PDEs reduction and $\mu$-symmery.} Note di matematica 23,n 2. 33-73, 2004.
\bibitem[Mu-Ro-2001]{[6]} {\sc Muriel, C.}, and {\sc Romero, J.L.} {\em New methods of reduction for ordinary differential equations}, IMA J. Appl. Math. 66, 111-125, 2001.
\bibitem[MuRo-2001]{[7]} {\sc Muriel, C.}, and {\sc Romero, J.L.} {\em $C^\infty$-Symmetries and non-solvable symmetry algebras}, IMA J. App. Math. 66,
441-498, 2001.
\bibitem[Mu-Ro-Ol-2006]{[14]}{\sc Muriel, C.}, and {\sc Romero, J.L.}, {\sc Olver,
P.}{\em Variational $c^\infty$-symmetries and Euler-Lagrange
equations}, Journal of differential equatiosn, 2006, 164-184.
\bibitem[Ol-1986]{[2]} {\sc Olver, P.J.} {\em Applications of Lie Groups to Differential Equations}, New York, Springer, 1986.
\bibitem[Ol-1995]{[3]} {\sc Olver, P.J.} {\em Equivalence, Invariants and Symmetry}, Cambridge University Press, 1995.
\bibitem[St-1989]{[4]} {\sc Stephani, H.} {\em Differential Equations}, Cambridge University Press, 1989.
\bibitem[Va-2003]{[9]} {\sc Van Brunt, B.} {\em The Calculus of variations}, Springer, 2003.
\end{thebibliography}
\end{document}